\documentclass{ws-ijnt}

\begin{document}


%
\catchline{}{}{}{}{}
%

\title{The $p$-torsion of curves with large $p$-rank}

\author{Rachel Pries}

\address{Mathematics, Colorado State University\\ 
Fort Collins, CO, 80523-1874, USA\\ 
\email{pries@math.colostate.edu}}

\maketitle

\begin{history}
\received{(Day Month Year)}
\accepted{(Day Month Year)}
\comby{xxx}
\end{history}

\begin{abstract}
Consider the moduli space of smooth curves of genus $g$ and $p$-rank $f$
defined over an algebraically closed field $k$ of characteristic $p$.
It is an open problem to classify which group schemes occur 
as the $p$-torsion of the Jacobians of these curves for $f < g-1$.
We prove that the generic point of every component of this moduli space has $a$-number 1 when $f=g-2$ and $f=g-3$.
Likewise, we show that a generic hyperelliptic curve with $p$-rank $g-2$ has $a$-number 1 when $p \geq 3$. 
We also show that the locus of curves with $p$-rank $g-2$ and $a$-number $2$ 
is non-empty with codimension $3$ in ${\cal M}_g$ when $p \geq 5$.
We include some other results when $f=g-3$.
The proofs are by induction on $g$ while fixing $g-f$.
They use computations about certain components of the boundary of ${\cal M}_g$.
\end{abstract}



\section{Introduction}	

Consider the moduli space ${\cal A}_g$ of
principally polarized abelian varieties of dimension $g$
which are defined over an algebraically closed field $k$ of characteristic $p$.
There are several important stratifications of ${\cal A}_g$ which characterize the $p$-torsion of the corresponding abelian varieties.
Very little is known about how these strata intersect the Torelli locus of Jacobians of curves.
It is an open question whether the Torelli locus intersects these strata and, if so, 
the extent to which the intersection is transversal.
In particular, one can ask which group schemes occur as the $p$-torsion $J_X[p]$ 
of the Jacobian $J_X$ of a curve $X$ of genus $g$.

The $p$-rank $f_X$ is the one invariant of $J_X[p]$ for which this question is moderately well-understood.  
Recall that $f_X={\rm dim}_{{\mathbb F}_p} {\rm Hom}(\mu_p, J_X[p])$ where the group scheme 
$\mu_p$ is the kernel of Frobenius on ${\mathbb G}_m$.
One can often compute the $p$-rank of a fixed curve $X$ for a fixed prime $p$, e.g., \cite{Y} when $X$ is hyperelliptic. 
It was recently proven that there exist smooth curves of genus $g$ 
with every possible $p$-rank \cite[Thm.\ 2.3]{FVdG}.
Likewise, there exist smooth hyperelliptic curves of genus $g$ with every possible $p$-rank, 
\cite[Thm.\ 1]{GP:hyper} for $p \geq 3$ and \cite[Thm.\ 3]{Z:hypeveryprank} for $p=2$. 

It is now natural to try to classify the group schemes which occur as the $p$-torsion for a smooth $k$-curve of
genus $g$ and $p$-rank $f$.  These are of the form $({\mathbb Z}/p \oplus \mu_p)^f \oplus {\mathbb G}$ 
where there are $2^{g-f-1}$ possibilities for ${\mathbb G}$ if $f \leq g-1$ \cite{Kraft, O:strat}.  
One can sometimes distinguish these by the $a$-number  
$a_X={\rm dim}_k {\rm Hom}(\alpha_p, J_X[p])$ where $\alpha_p$ is the kernel of Frobenius on ${\mathbb G}_a$.
Among these, there is a unique choice of ${\mathbb G}$ with $a$-number $1$, which we denote $I_{g-f}$.  
Precisely, $I_{g-f}$ is the unique symmetric ${\rm BT}_1$ group scheme of rank $p^{2(g-f)}$, $p$-rank 0, and $a$-number 1.
For example, $I_1$ is the $p$-torsion of a supersingular elliptic curve.
The group scheme $I_2$ occurs as the $p$-torsion of a supersingular non-superspecial abelian surface.

In this paper, we conjecture that the curve corresponding to the generic point of any 
component of the locus in ${\cal M}_g$ of smooth curves of genus $g$ with $p$-rank $f$ 
has group scheme $J_X[p] \cong ({\mathbb Z}/p \oplus \mu_p)^f \oplus I_{g-f}$ 
and thus $a$-number $1$ for $f \leq g-1$.  
It is easy to see that this conjecture is true when $f=g-1$.
When $f=g-2$ and $f=g-3$, we prove the conjecture along with some other results.
The results are stated more precisely in the paper in terms of the dimension of the
intersection of loci of ${\cal M}_g$ and ${\cal H}_g$ but they imply the following.
Suppose $g \geq 3$.

\begin{description}
\item{{\bf Theorem \ref{Cf2}:}} 
There is a family (with dimension $3g-5$) of smooth curves of genus $g$ with $p$-rank $g-2$ and $a$-number $1$.
\medskip
\item{{\bf Theorem \ref{Cf2hyper}:}}
If $p \geq 3$, there is a family (with dimension $2g-3$) of smooth hyperelliptic curves of genus $g$
with $p$-rank $g-2$ and $a$-number $1$.  
\medskip
\item{{\bf Theorem \ref{Cf3}:}}
There is a family (with dimension $3g-6$) of smooth curves of genus $g$ with $p$-rank $g-3$ and $a$-number $1$.
\medskip
\item{{\bf Corollary \ref{Ct2}:}} 
If $p \geq 5$, there is a family (with dimension $3g-6$) of smooth curves of genus $g$ with $p$-rank $g-2$ and $a$-number two.
\medskip
\item{{\bf Proposition \ref{PuseR}:}} For $p \geq 3$ and $g$ odd, $g \not \equiv 1 \bmod p$, and $g>6(p-1)$
all four possibilities of group scheme occur as $J_X[p]$ for a smooth curve $X$ of genus $g$ and $p$-rank $g-3$.
\end{description}

The proofs are by induction on $g$ with the quantity $r=g-f$ fixed.  
The induction result works for any value of $r$.  
It involves the computation of the dimension of certain components of the boundary of ${\cal M}_g$.
For the initial case of the induction process, one needs results about the locus of curves of genus 
$r$ with $p$-rank $0$.  At the moment such results are known only when $r \leq 3$. 
For example, we could prove the conjecture when $f=g-4$ with more information about 
smooth curves with genus $r = 4$ and $p$-rank $0$. 

Here is an outline of the paper. Section \ref{Sbackground} contains background.  
The inductive result is in Section \ref{Sboundary}.  The results for large $p$-rank are in Section \ref{Sbigf}.
In Section \ref{Sremark}, we explain the limitations of the inductive approach used in this paper 
and show that some of the hypotheses of the results are necessary.

One can interpret the results of this paper as 
evidence that the Torelli locus is in general position relative to the stratification of ${\cal A}_g$ by $p$-torsion type.
This positioning should not be taken for granted, e.g., \cite{OVdG}.  
The theory of tautological classes may play a fundamental role in investigating these questions further, \cite{V:cycles}.
Finally, we note that not all group schemes occur for small $p$, e.g., \cite{Re, Z:p2nohypsup}, 
or for curves with nontrivial automorphism group, e.g., \cite{B, El:bound}.

\section{Background} \label{Sbackground}

\subsection{Moduli spaces} \label{Smoduli}

Let $k$ be an algebraically closed field of characteristic $p >0$.
Let ${\cal M}_g$ (resp.\ ${\cal H}_g$) denote the moduli stack of smooth projective irreducible (resp.\ and hyperelliptic) $k$-curves 
of genus $g$.
Let ${\cal A}_g$ denote the moduli stack of principally polarized abelian varieties of dimension $g$ over $k$.
We assume $g \geq 2$ to avoid trivial cases.
Recall that ${\cal M}_g$ has dimension $3g-3$ and ${\cal H}_g$ has dimension $2g-1$.

These moduli stacks have natural compactifications $\overline{{\cal M}}_g$ (resp.\ $\overline{{\cal H}}_g$) whose points 
parametrize stable (resp.\ and hyperelliptic) curves of genus $g$.
Consider also the toroidal compactification $\overline{{\cal A}}_g$ of ${\cal A}_g$ \cite{FC}.
Then $\overline{{\cal M}}_g$, $\overline{{\cal H}}_g$, and $\overline{{\cal A}}_g$ are smooth proper algebraic stacks.
By \cite{V:stack}, a smooth stack has the same intersection-theoretic properties as a smooth scheme.
In particular, since $k$ is algebraically closed, we have the following intersection property:
if two closed substacks of $\overline{{\cal M}}_g$ (or $\overline{{\cal H}}_g$, or $\overline{{\cal A}}_g$) intersect 
then the codimension of their intersection is at most the sum of their codimensions.

Without further comment, a (generic) point means a (generic) geometric point.
We identify a curve with the point of ${\cal M}_g$ corresponding to its isomorphism class (and likewise for ${\cal H}_g$ and ${\cal A}_g$).

\subsection{The $p$-rank and $a$-number}

Let $X$ be a $k$-curve of genus $g$.  
The Jacobian $J_X$ of $X$ corresponds to a $k$-point of $\overline{{\cal A}}_g$.
The $p$-torsion $J_X[p]$ is a group scheme of rank $p^{2g}$. 
Recall that the group scheme $\mu_p$ is the kernel of Frobenius on ${\mathbb G}_m$
and the group scheme $\alpha_p$ is the kernel of Frobenius on ${\mathbb G}_a$.

Two invariants of the $p$-torsion of an abelian variety are the $p$-rank and $a$-number.
Let $f_X={\rm dim}_{{\mathbb F}_p} {\rm Hom}(\mu_p, J_X[p])$ be the {\it $p$-rank} of $X$.  Then $0 \leq f_X \leq g$.
Let $a_X={\rm dim}_k {\rm Hom}(\alpha_p, J_X[p])$ be the {\it $a$-number} of $X$.
Then $0 \leq a_X \leq g-f_X$.
The $p$-rank can only decrease under specialization, while the $a$-number can only increase.

Given $g$ and $f$ such that $0 \leq f \leq g$, let $V_{g,f}$ denote the closed substack of $\overline{{\cal M}}_g$ 
consisting of curves $X$ of genus $g$ with $f_X \leq f$. 
The locus $V_{g,f}$ is pure of dimension $2g-3+f$ by \cite[Thm.\ 2.3]{FVdG}. 
Also $V_{g,f} \cap {\cal H}_g$ is pure of dimension $g-1+f$ by \cite[Thm.\ 1]{GP:hyper} when $p \geq 3$
and by \cite[Thm.\ 3]{Z:hypeveryprank} when $p=2$.
In other words, $V_{g,f}$ has codimension $g-f$ in ${\cal M}_g$ and in ${\cal H}_g$.

In particular, the generic point $X$ of ${\cal M}_g$ or of ${\cal H}_g$ 
has $f_X=g$, and thus $a_X=0$ and $J_X[p] \cong ({\mathbb Z}/p \oplus \mu_p)^g$; a curve with these attributes is {\it ordinary}. 
The non-ordinary locus $V_{g,g-1}$ has codimension one in ${\cal M}_g$ and in ${\cal H}_g$.
The generic point $X$ of each component of $V_{g,g-1}$ and $V_{g,g-1} \cap {\cal H}_g$
has $f_X=g-1$, and thus $a_X=1$ and $J_X[p] \cong ({\mathbb Z}/p \oplus \mu_p)^{g-1} \oplus I_1$;
a curve with these attributes is {\it almost ordinary}.
Here $I_1$ is the $p$-torsion of a supersingular elliptic curve (see Example \ref{EI1}). 

Let $T_{g,a}$ denote the closed substack of $\overline{{\cal M}}_g$ consisting of curves $X$ of genus $g$ with $a_X \geq a$.
The results above show that $T_{g,0}$ has codimension $0$ in ${\cal M}_g$ (resp.\ ${\cal H}_g$)
and $T_{g,1}$ has codimension $1$ in ${\cal M}_g$ (resp.\ ${\cal H}_g$). 
There are not many other results about the intersection of $T_{g,a}$ with ${\cal M}_g$ or ${\cal H}_g$.
This paper includes a result for $a=2$.   

Finally, let $V'_{g,f}$ (resp.\ $T'_{g,2}$) denote the closed substack of $\overline{{\cal A}}_g$ 
consisting of abelian varieties with $p$-rank at most $f$ (resp.\ 
$a$-number at least two).

\subsection{Group schemes} \label{Sgroupscheme}

The results in this paper can be expressed in terms of the $p$-rank and $a$-number,
or in terms of the isomorphism class of the $p$-torsion, which is a group scheme.
These group schemes were classified by Kraft (unpublished) \cite{Kraft} and independently by Oort \cite{O:strat}.
Together with Ekedahl, Oort used this classification to stratify ${\cal A}_g$. 
Here is a brief explanation of this topic. 
A more complete historical and mathematical description can be found in \cite{M:group, O:strat, P:groupscheme}.

The group schemes ${\mathbb G}$ which occur as the $p$-torsion of a principally polarized abelian variety defined over $k$
are symmetric ${\rm BT}_1$ group schemes (short for truncated Barsotti-Tate group of level 1).
A ${\rm BT}_1$ is a finite commutative $k$-group scheme annihilated by $p$
having actions by semi-linear operators Frobenius $F$ and Verschiebung $V$ 
that satisfy certain properties, see \cite[2.1]{O:strat}. 

The isomorphism type of a symmetric ${\rm BT}_1$ group scheme ${\mathbb G}$ can be encapsulated into combinatorial data.
If ${\mathbb G}$ has rank $p^{2g}$, 
then there is a {\it final filtration} $N_1 \subset N_2 \subset \ldots \subset N_{2g}$ of the covariant Dieudonn\'e module of 
${\mathbb G}$ which is stable under $V$ and $F^{-1}$ with ${\rm rank}(N_i)=p^i$ \cite[2.2]{O:strat}.
The {\it Ekedahl-Oort type} of ${\mathbb G}$ is $[\nu_1, \ldots, \nu_g]$ where $p^{\nu_i}={\rm rank}(V(N_i))$.
There are some conditions on the filtration: e.g., $\nu_{i+1} \leq \nu_i +1$ \cite[Lemma 2.4]{O:strat}.
The $p$-rank is ${\rm max}\{i \ | \ \nu_i=i\}$ and the $a$-number is $g-\nu_g$.

Using this combinatorial data, 
one can see that there are $2^g$ isomorphism types of symmetric ${\rm BT}_1$ group schemes of rank $p^{2g}$
and that $2^{g-f-1}$ of these have $p$-rank $f$.    
There is a stratification of ${\cal A}_g$ by Ekedahl-Oort type. 
The dimension of the stratum of ${\cal A}_g$ whose points have Ekedahl-Oort type $[\nu_1, \ldots, \nu_g]$ is $\sum_{i=1}^g \nu_i$
\cite[Thm. 1.2]{O:strat}.

\subsection{The boundary of ${\cal M}_g$ and ${\cal H}_g$}

The boundary $\overline{{\cal M}}_g -{\cal M}_g$ consists of components
$\Delta_0$ and $\Delta_i$ for integers $1 \leq i \leq g/2$.
Each boundary component has codimension 1 in $\overline{{\cal M}}_g$.
The boundary component $\Delta_0$ will be useful in the next section.

The generic point of $\Delta_0$ corresponds to the isomorphism class of an irreducible curve 
$X_0'$ with one ordinary double point $P'$. 
The normalization of $X_0'$ is a smooth curve $X_0$ of genus $g-1$, 
and the inverse image of $P'$ consists of two distinct points $P_1,P_2$ of $X_0$.
Recall that $J_{X_0'}$ is a semi-abelian variety which is an extension of $J_{X_0}$ by a torus, \cite[pg. 246]{BLR}.
This torus contains a copy of $\mu_p$.  It follows that $f_{X_0'} =f_{X_0} +1$.  Also $a_{X_0'}=a_{X_0}$.

Likewise, if $p \geq 3$,  the boundary $\overline{{\cal H}}_g -{\cal H}_g$ consists of components
$\Delta_0$ and $\Delta_i$ for $1 \leq i \leq g/2$.
Each boundary component has codimension 1 in $\overline{{\cal H}}_g$.
A {\it singular hyperelliptic curve} is a curve $X$ corresponding to a point of $\overline{{\cal H}}_g -{\cal H}_g$. 
If $X_0' \in \Delta_0$ is hyperelliptic, then so is its normalization $X_0$ and the two points $P_1$ and $P_2$ 
are exchanged by the hyperelliptic involution.


A complete substack of $\tilde{{\cal M}}_g=\overline{{\cal M}}_g - \Delta_0$ 
has codimension at least $g$ \cite[pg 80]{D:complete} \cite{L:complete}, see also \cite[2.4]{FVdG}.
For $p \geq 3$, a complete substack of $\tilde{{\cal H}}_g=\overline{{\cal H}}_g -\Delta_0$ 
has codimension at least $g$ \cite[2.6]{FVdG}
In other words, if $L$ is a closed substack of codimension less than $g$ in either $\overline{{\cal M}}_g$ or $\overline{{\cal H}}_g$
then $L$ intersects $\Delta_0$.

\section{Conjecture and Inductive Result} \label{Sboundary}

\subsection{Description of some group schemes} \label{Sgrouplemma}

The group scheme $I_r$ found in the next lemma plays an important role in this paper.

\begin{lemma} \label{Lunique}
For $r \in {\mathbb N}^+$, there is a unique symmetric ${\rm BT}_1$ group scheme with rank $p^{2r}$, $p$-rank 0, and $a$-number 1. 
It has Ekedahl-Oort type $[0,1, \ldots, r-1]$.
We denote this group scheme by $I_r$.
\end{lemma}

\begin{proof}
This is proved in \cite[Lemma 3.1]{P:groupscheme}.  For the convenience of the reader, here are a few details of the proof.
Let $I_r$ be a symmetric ${\rm BT}_1$ group scheme with rank $p^{2r}$, $p$-rank 0, and $a$-number 1.
It is sufficient to show that the Ekedahl-Oort type of $I_r$ is uniquely determined.

The $p$-rank 0 condition implies that $V$ acts nilpotently on $I_r$, so $\nu_1=0$.
The $a$-number 1 condition implies that $r-1$ is the dimension of $V^2 I_r$, so $\nu_r=r-1$.
The numerical conditions $\nu_{i+1} \leq \nu_i +1$ from \cite[Lemma 2.4]{O:strat} imply 
that there is a unique possibility for the Ekedahl-Oort type of $I_r$, namely $[0,1, \ldots, r-1]$.
\end{proof}

\begin{example} \label{EI1}
The group scheme $I_1$ occurs as the $p$-torsion of a supersingular elliptic curve.
By \cite[Ex.\ A.3.14]{G:book}, $I_1$ fits into a non-split exact sequence of the form $0 \to \alpha_p \to I_1 \to \alpha_p \to 0$.
The image of the embedded $\alpha_p$ is unique and is the kernel of both Frobenius and Verschiebung.  
\end{example}
 
\begin{example} \label{EI2}
The group scheme $I_2$ occurs as the $p$-torsion of a supersingular non-superspecial abelian surface. 
By \cite[Ex.\ A.3.15]{G:book}, there is a filtration
$H_1 \subset H_2 \subset I_2$ where $H_1 \cong \alpha_p$,
$H_2/H_1 \cong \alpha_p \oplus \alpha_p$, and $I_2/H_2 \cong \alpha_p$.
There is an exact sequence $0 \to H_1 \to G_1 \oplus G_2 \to H_2 \to 0$ where
$G_1$ (resp.\ $G_2$) is the kernel of Frobenius (resp.\ Verschiebung).  Also $G_1 \subset H_2$ and $G_2 \subset H_2$.
The $p$-rank of $I_2$ is zero since ${\rm Hom}(\mu_p, I_2)=0$.
The $a$-number of $I_2$ is one since ${\rm ker}(V^2)=G_1 \oplus G_2$ has rank $p^3$.
\end{example}

\begin{remark}  
One can also describe group schemes using the theory of covariant Dieudonn\'e modules.
Briefly, let $\sigma$ denote the Frobenius automorphism of $k$. 
Consider the non-commutative ring $E=k[F,V]$ with the relations 
$FV=VF=p$ and $F\lambda=\lambda^\sigma F$ and $\lambda V=V \lambda^\sigma$ for all $\lambda \in k$.
There is an equivalence of categories between ${\rm BT}_1$ group schemes 
${\mathbb G}$ and finite left $E$-modules $D({\mathbb G})$
The covariant Dieudonn\'e module for the group scheme $I_r$ is $E/E(F^r-V^r)$ \cite[Lemma 3.1]{P:groupscheme}.
\end{remark}

Here are a few facts about the Ekedahl-Oort strata of ${\cal A}_3$.
First, there are four possibilities for the isomorphism class of 
a symmetric ${\rm BT}_1$ group scheme ${\mathbb G}'$ with rank $p^6$ and $p$-rank $0$: 
(i) ${\mathbb G}'=I_3$; (ii) ${\mathbb G}'= I_{3,2}$; (iii) ${\mathbb G}'=I_2 \oplus I_1$; and (iv) ${\mathbb G}'=(I_1)^3$. 
Here $I_{3,2}$ is the symmetric ${\rm BT}_1$ group scheme of rank $p^6$, $p$-rank $0$, and $a$-number $2$ 
which is not $I_2 \oplus I_1$.
One can show that $I_{3,2}$ has Ekedahl-Oort type $[0,1,1]$ and covariant Dieudonn\'e module $E/E(F-V^2) \oplus E/E(V-F^2)$ 
\cite[Lemma 3.4]{P:groupscheme}.
The group scheme $I_3$ is the only one of the four with $a$-number 1.  

Secondly, here is a useful lemma about the $a$-number $2$ locus of ${\cal A}_3$.  

\begin{lemma} \label{Ldim3a2}
In ${\cal A}_3$, every component of $T'_{3,2}$ has dimension $3$ and 
every component of $V'_{3,0} \cap T'_{3,2}$ has dimension $2$.  
\end{lemma}

\begin{proof}
Consider the Ekedahl-Oort types occuring in dimension three with $a$-number at least $2$, namely 
$\nu_1=[1,1,1]$, $\nu_2=[0,1,1]$, $\nu_3=[0,0,1]$, and $\nu_4=[0,0,0]$.  
For each of these types $\nu$, let $L_\nu$ be the locus of ${\cal A}_3$ of abelian threefolds with Ekedahl-Oort type $\nu$.
By \cite[Thm.\ 1.2]{O:strat}, each of these loci $L_\nu$ is non-empty and quasi-affine; all components of $L_{\nu}$ have the same 
dimension, which is respectively $3$, $2$, $1$, and $0$.  
The only types that occur under the additional restriction of having $p$-rank $0$ are $\nu_2$, $\nu_3$, $\nu_4$.  

By \cite[Prop.\ 11.1]{O:strat}, if $1 \leq i \leq 3$, then $L_{\nu_{i+1}}$ is in the closure of $L_{\nu_i}$ in ${\cal A}_3$.
Thus the generic point of every component of $T'_{3,2}$ is contained in $L_{\nu_1}$ and 
the generic point of every component of $V'_{3,0} \cap T'_{3,2}$ is contained in $L_{\nu_2}$.
Thus every component of $T'_{3,2}$ has dimension $3$ and 
every component of $V'_{3,0} \cap T'_{3,2}$ has dimension $2$.  
\end{proof}

\subsection{Conjecture}

\begin{conjecture} \label{Cgeneric}
Let $0 \leq f < g$.
Let $X$ be the generic point of any component of $V_{g,f} \cap {\cal M}_g$ (resp.\ $V_{g,f} \cap {\cal H}_g$ for $p \geq 3$).
Then $J_X[p] \cong({\mathbb Z}/p \oplus \mu_p)^f \oplus I_{g-f}$.
Equivalently, a generic smooth (resp.\ and, for $p \geq 3$,  hyperelliptic) curve of genus $g$ and $p$-rank $f$ has $a$-number 1.
\end{conjecture}

The motivation for Conjecture \ref{Cgeneric} is that a generic point of $V'_{g,f}$ has $a$-number $1$ by \cite{O:strat}.
If Conjecture \ref{Cgeneric} is true for $V_{g,f} \cap {\cal H}_g$ 
then it is true for at least one component of $V_{g,f} \cap {\cal M}_g$.
The hypothesis $p \geq 3$ is necessary for the hyperelliptic case, see Section \ref{Sp=2}.

\subsection{Inductive result} \label{Sinduct}

\begin{proposition} \label{Pinduction} 
Suppose $g \geq 3$ and $1 \leq f < g$.
Suppose $g'$ is an integer such that $g-f  \leq g ' < g$.  Let $f'=g'-(g-f)$. 
\begin{description}
\item{(1)} If Conjecture \ref{Cgeneric} is true for $V_{g',f'} \cap {\cal M}_{g'}$ then it is true for $V_{g,f} \cap {\cal M}_{g}$.
\item{(2)} If $p \geq 3$ and if Conjecture \ref{Cgeneric} is true for $V_{g',f'} \cap {\cal H}_{g'}$ 
then it is true for $V_{g,f} \cap {\cal H}_{g}$.
\end{description}
\end{proposition}

\begin{proof}
The proof is by induction on $g$ with the quantity $r=g-f$ fixed.
Note that $g-f=g'-f'$. 
Without loss of generality, one can take the inductive hypothesis to be that Conjecture \ref{Cgeneric}
is true when $g'=g-1$ and $f'=f-1$. 

\begin{description}
\item{(1)}
Let $X$ be the generic point of a component $\Gamma$ of $V_{g,f} \cap {\cal M}_g$.
Now, ${\rm dim}(\Gamma)=2g-3+f$ and
${\rm dim}(V_{g, f-1} \cap {\cal M}_g)=2g-4+f$ by \cite[Thm 2.3]{FVdG}.  
Thus $f_X=f$.

If $a_X \geq 2$, then $\Gamma \subset T_{g,2}$.
Let $\overline{\Gamma}$ be the closure of $\Gamma$ in $\overline{{\cal M}}_g$.
Then ${\rm codim}(\overline{\Gamma}, \overline{{\cal M}}_g) < g$ since $f \geq 1$.
By \cite[pg 80]{D:complete} and \cite{L:complete}, $\overline{\Gamma}$ intersects $\Delta_0$. 
Since $\Delta_0$ has codimension $1$ in $\overline{{\cal M}}_g$, by Section \ref{Smoduli},
\begin{equation*}{\rm dim}(\overline{\Gamma}) \leq {\rm dim}(V_{g,f} \cap T_{g,2} \cap \Delta_0)+1.\end{equation*}  
A generic point of $V_{g,f} \cap T_{g,2} \cap \Delta_0$ corresponds to a smooth curve $X_0$ of genus $g-1$ 
with $2$ points identified in an ordinary double point.  Also $f_{X_0}=f-1$ and $a_{X_0} \geq 2$.  
Thus, 
\begin{equation*}{\rm dim}(V_{g,f} \cap T_{g,2} \cap \Delta_0)={\rm dim}(V_{g-1,f-1} \cap T_{g-1,2} \cap \overline{{\cal M}}_{g-1}) +2. \end{equation*}
By the inductive hypothesis, the generic point of each component of $V_{g-1,f-1} \cap \overline{{\cal M}}_{g-1}$ has $a$-number 1.
It follows that 
\begin{equation*}{\rm dim}(V_{g-1,f-1} \cap T_{g-1,2} \cap \overline{{\cal M}}_{g-1}) < 2(g-1)-3+(f-1).\end{equation*}
Then ${\rm dim}(\overline{\Gamma}) < 2g+f-3$ which is a contradiction.  
Thus $a_X=1$ and $J_X[p] \cong({\mathbb Z}/p \oplus \mu_p)^f \oplus I_{g-f}$ by Lemma \ref{Lunique}.

\item{(2)}  The proof is similar to part (1).
Let $X$ be the generic point of a component $\Gamma$ of $V_{g,f} \cap \overline{{\cal H}}_g$.
Then ${\rm dim}(\Gamma)=g-1+f$ and $f_X=f$ by \cite[Thm.\ 1]{GP:hyper}.
By \cite[Lemma 2.6]{FVdG}, $\overline{\Gamma}$ intersects $\Delta_0$ when $p \geq 3$.
If $a_X \geq 2$, then
\begin{equation*} {\rm dim}(\Gamma) \leq {\rm dim}(\overline{{\cal H}}_g \cap V_{g,f} \cap T_{g,2} \cap \Delta_0)+1.\end{equation*}
A generic point of the latter corresponds to a smooth hyperelliptic curve $X_0$ of genus $g-1$ 
with $2$ points identified.  These two points form one orbit under the hyperelliptic involution.
Also $f_{X_0}=f-1$ and $a_{X_0} \geq 2$.  
Thus, 
\begin{equation*}{\rm dim}(\overline{{\cal H}}_g \cap V_{g,f} \cap T_{g,2} \cap \Delta_0)
={\rm dim}(\overline{{\cal H}}_{g-1} \cap V_{g-1,f-1} \cap T_{g-1,2}) +1.\end{equation*}
By the inductive hypothesis, ${\rm dim}(\overline{{\cal H}}_{g-1} \cap V_{g-1,f-1} \cap T_{g-1,2}) < g-3 +f$, which yields a contradiction.
\end{description}
\end{proof}

\section{Results for the Case of Large $p$-Rank} \label{Sbigf}

This section includes results about the group schemes which occur
as the $p$-torsion of Jacobians of smooth curves of genus $g$ and $p$-rank $f$ when $f$ is large relative to $g$.
For $p$-rank $f=g$, the only group scheme possible is $({\mathbb Z}/p \oplus \mu_p)^g$
and this occurs for the generic point of ${\cal M}_g$ and ${\cal H}_g$.
For $p$-rank $f=g-1$, the only group scheme possible is $({\mathbb Z}/p \oplus \mu_p)^{g-1} \oplus I_1$
and this occurs with codimension one in ${\cal M}_g$ and ${\cal H}_g$.

If $r=g-f \geq 1$, there are $2^{r-1}$ possibilities for the group scheme $J_X[p]$ of a smooth curve with $p$-rank $f$.
By the inductive process in Section \ref{Sboundary}, the existence of curves of genus $r$ and $p$-rank $0$ can 
yield results about curves of genus $g$ with $p$-rank $g-r$.  
The situation when $r \leq 3$ is well understood due to the fact that ${\rm dim}({\cal M}_r)={\rm dim}({\cal A}_r)$ in these cases.
We amalgamate these ideas to deduce results when $f=g-2$ and $f=g-3$.

As $f$ decreases relative to $g$, the situation becomes more difficult to study.  This is true first because there  
are an increased number of possibilities for $J_X[p]$ and second because of the increase in codimension of the 
relevant loci in ${\cal M}_g$, see Section \ref{Scodim}.

\subsection{The conjecture is true for ${\cal M}_g$ and ${\cal H}_g$ when $f=g-2$} \label{Sf2}

Suppose $g \geq 2$ and $f=g-2$.
Then $J_X[p]$ is isomorphic to either $({\mathbb Z}/p \oplus \mu_p)^{g-2} \oplus I_2$ or 
$({\mathbb Z}/p \oplus \mu_p)^{g-2} \oplus (I_1)^2$, which have $a$-number 1 and 2 respectively.  
The next results imply that the locus of curves with $p$-rank $g-2$ and $a$-number 1 
is non-empty with codimension $2$ in ${\cal M}_g$ and in ${\cal H}_g$.     

\begin{theorem} \label{Cf2}
The group scheme $({\mathbb Z}/p \oplus \mu_p)^{g-2} \oplus I_2$ occurs as the $p$-torsion of 
the generic point of every component of $V_{g,g-2} \cap {\cal M}_g$.
So, there is a family of dimension $3g-5$ consisting of smooth curves 
of genus $g$ with $p$-rank $g-2$ and $a$-number $1$.
\end{theorem}

\begin{proof}
Every component of $V_{2,0} \cap {\cal M}_2$ has dimension $1$.  
There are only finitely many points of ${\cal A}_2$ with $a$-number 2 \cite[Thm.\ 2.10]{O:sup2}.
Thus the generic point of each component of $V_{2,0} \cap {\cal M}_2$ has $a$-number $1$, 
and group scheme $I_2$ by Lemma \ref{Lunique}.  
This means that Conjecture \ref{Cgeneric} is true for $V_{2,0} \cap {\cal M}_2$.
The result for $V_{g,g-2} \cap {\cal M}_g$ then follows from Proposition \ref{Pinduction}(1).
\end{proof}

\begin{theorem} \label{Cf2hyper}
Let $p \geq 3$.  The group scheme $({\mathbb Z}/p \oplus \mu_p)^{g-2} \oplus I_2$ occurs as the $p$-torsion of 
the generic point of every component of $V_{g,g-2} \cap {\cal H}_g$.
So, there is a family of dimension $2g-3$ consisting of smooth hyperelliptic curves 
of genus $g$ with $p$-rank $g-2$ and $a$-number $1$.
\end{theorem}

\begin{proof}
The proof is the same as for Theorem \ref{Cf2} and uses Proposition \ref{Pinduction}(2).
\end{proof}

\subsection{The conjecture is true for ${\cal M}_g$ when $f=g-3$} \label{Sf3}

Suppose $g \geq 3$ and $f=g-3$.
The next result is that the locus of curves with $p$-rank $g-3$ and $a$-number 1 
is non-empty with codimension $3$ in ${\cal M}_g$.

\begin{theorem} \label{Cf3}
The group scheme $({\mathbb Z}/p \oplus \mu_p)^{g-3} \oplus I_3$ occurs as the $p$-torsion of 
the generic point of every component of $V_{g,g-3} \cap {\cal M}_g$.
So, there is a family of dimension $3g-6$ consisting of smooth curves 
of genus $g$ with $p$-rank $g-3$ and $a$-number $1$.
\end{theorem}

\begin{proof}
Every component of $V_{3,0} \cap {\cal M}_g$ has dimension 3 and its image in ${\cal A}_3$ under the Torelli map also has dimension 3.
On the other hand, ${\rm dim}(V'_{3,0} \cap T'_{3,2}) = 2$ by Lemma \ref{Ldim3a2}.  
Thus the generic point of each component of $V_{3,0} \cap {\cal M}_g$ has $a$-number $1$, and group scheme 
$I_3$ by Lemma \ref{Lunique}.  
Thus Conjecture \ref{Cgeneric} is true for $V_{3,0} \cap {\cal M}_3$.
The result follows from Proposition \ref{Pinduction}(1).
\end{proof}

\begin{remark}
In \cite[Cor.\ 3.2]{NR:nonhyper}, for any finite field $k$ of characteristic $2$, the authors calculate the number of 
$k$-isomorphism classes of curves of genus $3$ with a given Newton polygon.
From the data one sees that $V_{3,0} \cap {\cal M}_3$ is irreducible when $p=2$.
\end{remark}

\subsection{The case of $p$-rank $g-2$ and $a$-number $2$} \label{Sf2a2} 

Suppose $g \geq 3$ and $f=g-2$.
In this section, we consider the group scheme $({\mathbb Z}/p \oplus \mu_p)^{g-2} \oplus (I_1)^2$.  
The next result shows that the locus of curves with this group scheme has codimension $3$ in $\overline{{\cal M}}_g$ for $g \geq 3$.

\begin{corollary} \label{Ct2} 
Let $g \geq 3$.
\begin{description}
\item{(1)}  The group scheme $({\mathbb Z}/p \oplus \mu_p)^{g-2} \oplus (I_1)^2$ occurs as the $p$-torsion of 
the generic point of every component of $T_{g,2} \cap \overline{{\cal M}}_g$.
\item{(2)} If $p \geq 5$, there is a family of dimension $3g-6$
consisting of smooth curves of genus $g$ with $p$-rank $g-2$ and $a$-number $2$.
\end{description}
\end{corollary}

\begin{proof}
\begin{description}
\item{(1)}
Let $\overline{\Gamma}$ be a component of $T_{g,2}$.
The generic point $X$ of $\overline{\Gamma}$ has $2 \leq a_X \leq g-f_X$.
Since $({\mathbb Z}/p \oplus \mu_p)^{g-2} \oplus (I_1)^2$ is the unique symmetric ${\rm BT}_1$ group scheme of rank 
$p^{2g}$, $p$-rank $g-2$, and $a$-number 2, 
it suffices to show that $f_X=g-2$. 

The result is true when $g=3$.  To see this, recall from Lemma \ref{Ldim3a2} 
that every component of $T'_{3,2}$ has dimension $3$ with generic point having $p$-rank $1$ and $a$-number $2$.
The same is true for $\overline{\Gamma}$ since its image under the Torelli morphism is a component of $T'_{g,2}$.

The proof is now by induction on $g$.  Let $g \geq 4$.
Now ${\rm codim}(T'_{g,2}, \overline{{\cal A}}_g)=3$. 
Consider the intersection of $T'_{g,2}$ with the Torelli locus $\tau$ in $\overline{{\cal A}}_g$.
By Section \ref{Smoduli}, 
${\rm codim}(T'_{g,2} \cap \tau, \overline{{\cal A}}_g) \leq {\rm codim}(\tau, \overline{{\cal A}}_g) +3$
which yields ${\rm codim}(\overline{\Gamma}, \overline{{\cal M}}_g) \leq 3$.
Then ${\rm codim}(\overline{\Gamma}, \overline{{\cal M}}_g) =3$ by Theorem \ref{Cf2}. 

Since $3 < g$, \cite[Lemma 2.6]{FVdG} implies $\overline{\Gamma}$ intersects $\Delta_0$.
If $J_X[p] \not \cong ({\mathbb Z}/p \oplus \mu_p)^{g-2} \oplus (I_1)^2$
then $f_X \leq g-3$ and $X \in V_{g,g-3} \cap T_{g,2}$.
One checks that: 
\begin{align*}
{\rm dim}(V_{g,g-3} \cap \overline{\Gamma}) &  \leq {\rm dim}(V_{g,g-3} \cap T_{g,2} \cap \Delta_0) +1 \ {\rm and} \\
{\rm dim}(V_{g,g-3} \cap T_{g,2} \cap \Delta_0) & = {\rm dim}(V_{g-1, g-4} \cap T_{g-1,2} \cap \overline{{\cal M}}_{g-1})+2.
\end{align*}
By Theorem \ref{Cf3}, the generic point of every component of $V_{g-1, g-4} \cap \overline{{\cal M}}_{g-1}$ has $a$-number $1$.  Thus 
\begin{equation*}{\rm dim}(V_{g-1, g-4} \cap T_{g-1,2} \cap \overline{{\cal M}}_{g-1}) < {\rm dim}(V_{g-1,g-4} \cap \overline{{\cal M}}_{g-1})= 3g-9.\end{equation*}
Then ${\rm dim}(V_{g,g-3} \cap \overline{\Gamma}) < 3g-6$ which is too small to be generic in $\overline{\Gamma}$.
Thus $X$ must have group scheme $({\mathbb Z}/p \oplus \mu_p)^{g-2} \oplus (I_1)^2$.  
\item{(2)}
If $p \geq 5$, there exists a smooth curve $X$ of genus $g$ with $p$-rank $g-2$ and $a$-number $2$ by \cite[Cor. 4]{GP:hyper}.
The result follows since there is at least one component of $T_{g,2}$ which is not contained in the boundary of ${\cal M}_g$.
\end{description}
\end{proof}

\begin{remark}
There is an explicit formula for the number of smooth curves $X$ of genus $2$ with $f_X=0$ and $a_X=2$;
these curves exist for all $p \geq 5$ \cite[Section 3.3]{O:sup2}. 
\end{remark}

\begin{remark}
For each $i$ such that $1 \leq i \leq g/2$, 
there is a component of $T_{g,2}$, having dimension $3g-6$, 
which is contained in the boundary component $\Delta_i$ of $\overline{{\cal M}}_g$.  
\end{remark}

\subsection{The case of $p$-rank $g-3$ and $a$-number at least $2$} \label{Sf3a2} 

Suppose $g \geq 3$ and $f=g-3$.
Recall that if $X \in V_{g,g-3}$, then $J_X[p] \cong({\mathbb Z}/p \oplus \mu_p)^{g-3} \oplus {\mathbb G}'$ where there are four possibilities for ${\mathbb G}'$:
(i) ${\mathbb G}'=I_3$; (ii) ${\mathbb G}'= I_{3,2}$; (iii) ${\mathbb G}'=I_2 \oplus I_1$; and (iv) ${\mathbb G}'=(I_1)^3$. 
By Section \ref{Sf3}, case (i) occurs for the generic point of every component of $V_{g, g-3} \cap {\cal M}_g$.
This section includes partial results on the other three cases.

\begin{lemma} \label{Lg3other}
\begin{description}
\item{(ii)} For $p \geq 2$, there exists $Y \in {\cal M}_3$ with $J_Y[p]=I_{3,2}$;
\item({iii)} For $p \geq 3$, there exists $Y \in {\cal M}_3$ with $J_Y[p]=I_2 \oplus I_1$; 
\item{(iv)} For $p \geq 3$, there exists $Y \in {\cal M}_3$ with $J_Y[p]=(I_1)^3$.
\end{description}
\end{lemma}

\begin{proof}
\begin{description}
\item{(ii)} The locus of points of ${\cal A}_3$ with $p$-torsion $I_{3,2}$ has dimension two.  
Let $A$ be an abelian threefold with $A[p]=I_{3,2}$.
By the Torelli theorem, $A=J_Y$ for some curve $Y$ of genus $3$ with $p$-rank $0$.
Now $Y \not \in \Delta_0$ since $J_Y$ is an abelian variety with no toric part.
Since $I_{3,2}$ is indecomposable as a symmetric ${\rm BT}_1$ group scheme 
\cite[Lemma 3.4]{P:groupscheme}, 
$J_Y$ cannot be the product of abelian varieties of smaller dimension.
Thus $Y \not \in \Delta_1$ by \cite[pg.\ 246]{BLR}. 
Thus $Y$ is smooth and so $Y \in {\cal M}_3$.

\item{(iii)} By part (iv), there exists $Y' \in {\cal M}_3$ with $J_{Y'}[p]=(I_1)^3$.
The Torelli locus is open in ${\cal A}_3$ and contains $J_{Y'}$ by definition.
By \cite[Prop.\ 7.3]{O:strat}, there exists an irreducible dimension one sublocus $L \subset {\cal A}_3$ such that $J_{Y'} \in L$ 
and such that the generic point $\eta$ of $L$ has group scheme $I_2 \oplus I_1$.  
It follows that $\eta$ is contained in the Torelli locus.
Thus there exists $Y \in {\cal M}_3$ with $J_{Y}[p]=I_2 \oplus I_1$.  

\item{(iv)} This is proved in \cite[Theorem 5.12(1)]{O:hypsup}.
\end{description}
\end{proof}

\begin{proposition}  \label{PuseR}
Let $p \geq 3$.
Suppose $g$ is odd, $g \not \equiv 1 \bmod p$, and $g > 6(p-1)$.
Then there exists a curve $X \in V_{g,g-3} \cap {\cal M}_g$ 
so that $J_X[p] \cong ({\mathbb Z}/p \oplus \mu_p)^{g-3} \oplus {\mathbb G}'$ 
in each of the three cases:
(ii) ${\mathbb G}'= I_{3,2}$; (iii) ${\mathbb G}'=I_2 \oplus I_1$; and (iv) ${\mathbb G}'=(I_1)^3$.
\end{proposition}  

\begin{proof}
By Lemma \ref{Lg3other}, there exists a curve $Y \in {\cal M}_3$ with $J_Y[p]={\mathbb G}'$
if (ii) ${\mathbb G}'= I_{3,2}$; (iii) ${\mathbb G}'=I_2 \oplus I_1$; or (iv) ${\mathbb G}'=(I_1)^3$.
In each case, $f_Y=0$.

Let $\ell=(g-1)/2$.  By hypothesis, $p \nmid \ell$ and $\ell \geq 3(p-1)$.    
By \cite[Thm.\ 4.3.1]{R:etale}, under these two conditions on $\ell$, 
there exists an \'etale cyclic ${\mathbb Z}/\ell$-cover $X \to Y$ so that $f_X-f_Y=g_X-g_Y$.
By the Riemann-Hurwitz formula, $g_X= 1+ 2 \ell=g$.  Thus $f_X=g_X-3$.
Finally, $X$ dominates $Y$ so $J_Y[p] \subset J_X[p]$, completing the proof. 
\end{proof}

\begin{proposition}
Let $p=2$.  Suppose $g \equiv 1 \bmod 4$ and $g \geq 7$.  Then there exists a curve $X \in V_{g,g-3} \cap {\cal M}_g$ 
so that $J_X[p] \cong ({\mathbb Z}/p \oplus \mu_p)^{g-3} \oplus I_{3,2}$.
\end{proposition}

\begin{proof}
The proof is almost identical to that of Proposition \ref{PuseR}.
\end{proof}

\section{Questions and related results} \label{Sremark}

\subsection{Limitations of the boundary approach} \label{Scodim}

The inductive process of Section \ref{Sinduct} allows one to deduce results about $V_{g,f} \cap {\cal M}_{g}$ from results 
about $V_{r,0} \cap {\cal M}_{r}$ where $r=g-f$.  This produced results when $r \leq 3$.
Unfortunately, not much is known about $V_{r,0} \cap {\cal M}_{r}$ if $r \geq 4$.
One might ask whether the boundary of $V_{r,0} \cap {\cal M}_{r}$ could be used to deduce more information.
Here are some reasons why this is not straightforward.

First, $V_{r,0}$ does not intersect $\Delta_0$.
This is because the Jacobian of every curve in $\Delta_0$ has a toric part and so has $p$-rank at least $1$.

Second, the intersection of $V_{r,0}$ with the boundary component $\Delta_i$ is not useful for $1 \leq i \leq g/2$.
The generic point of $\Delta_i$ corresponds to the isomorphism class
of a singular curve $X$ with two irreducible components $X_i$ and
$X_{g-i}$, of genus $i$ and $g-i$, intersecting in an ordinary double point. 
Then $J_X[p] \cong J_{X_i}[p] \oplus J_{X_{g-i}}[p]$ \cite[pg.\ 246]{BLR}.
If $X \in V_{r,0} \cap \Delta_i$, then $f_{X_1}=f_{x_2}=0$ and $a_X=a_{X_1}+a_{X_2} \geq 2$.
In other words, there are no singular curves with $p$-rank $0$ and $a$-number $1$, which 
makes it difficult to determine whether a component of $V_{r,0}$ has $a$-number $1$.

One might also ask whether boundary methods could be used 
to deduce results about strata of ${\cal M}_g$ other than the $p$-rank strata.
For example, Corollary \ref{Ct2} is about the strata $T_{g,2} \cap {\cal M}_g$.
This becomes progressively more difficult 
for group schemes of rank $p^{2g}$ and $p$-rank $f$ as $g-f$ increases.
The first reason is that the number of possibilities for these group schemes increases exponentially in the quantity $g-f$.

The true difficulty, however, is geometric;
for some of these group schemes ${\mathbb G}$, the locus $L_{{\mathbb G}}$ of curves $X$ with $J_X[p] \cong {\mathbb G}$
is expected to have large codimension in $\overline{{\cal M}}_g$.
However, to guarantee that the closure of this locus intersects the boundary component $\Delta_0$
of ${\cal M}_g$, one requires that the codimension be less than $g$.
The base cases of the induction process thus involve curves of genus even larger than $r=g-f$.

For example, consider the group scheme ${\mathbb G}=({\mathbb Z}/p \oplus \mu_p)^{g-3} \oplus (I_1)^3$.
The locus $L_{{\mathbb G}}$ is expected to have codimension $6$ in $\overline{{\cal M}}_g$.
Thus, to produce curves with $f=g-3$ having group scheme $(I_1)^3$, one would need to 
start with results for curves of genus $6$.  
This was avoided in Section \ref{Sf3a2} by using \cite{R:etale}.

\subsection{Hyperelliptic curves of genus three} 

In this section, let $p \geq 3$ and consider the moduli space $V_{3,0} \cap {\cal H}_3$ of
hyperelliptic curves of genus $g$ and $p$-rank $0$.  
By \cite[Thm.\ 5.12(4)]{O:hypsup}, there exists $Y \in {\cal H}_3$ which is supersingular with $a_Y=1$ (and $f_Y=0$).
Then $J_Y[p] \cong I_3$ by Lemma \ref{Lunique}.
This implies that the generic point of {\it at least one} component of $V_{3,0} \cap {\cal H}_g$ has group scheme $I_3$. 
However, the inductive argument in this paper requires information about {\it every} component.

Here are some open questions whose solution would yield further progress.

\begin{question} \label{Qirred}
Is $V_{3,0} \cap {\cal H}_3$ irreducible?
\end{question}

If $p=3$, then $V_{3,0} \cap {\cal H}_3$ is irreducible by \cite[Prop.\ 3.5]{EP:g3small}.
Using this, the authors prove Conjecture \ref{Cgeneric} for $V_{g,g-3} \cap {\cal H}_g$ when $g \geq 3$ and $p=3$.
There do not seem to be any results about the number of components of $V_{3,0} \cap {\cal H}_3$ when $p \geq 5$.

\begin{question}
How does $T_{3,2} \cap {\cal H}_3$ intersect the supersingular locus of ${\cal H}_3$?
\end{question}

This relates to Conjecture \ref{Cgeneric} because every component of $V_{3,0} \cap {\cal H}_3$ 
has Newton polygon $G_{1,2} + G_{2,1}$ (three slopes of 1/3 and three slopes of 2/3) \cite[Thm.\ 1.12]{O:hypsup}.  
However, this Newton polygon cannot distinguish between the group schemes $I_3$
and $I_{3,2}$, which have $a$-numbers $1$ and $2$ respectively.
Here is some further evidence of the difficulty in comparing Newton polygons with $p$-torsion group schemes: 
if $J_X[p] \cong I_3$, then the Newton polygon of $X$ is typically $G_{1,2}+G_{2,1}$ 
but by \cite[Thm.\ 5.12 (2,4)]{O:hypsup} it can also be $3G_{1,1}$ (supersingular).

\begin{question} \label{QT2}
If $p \geq 5$, is the $p$-rank of the generic point of each two-dimensional component of ${\cal H}_3 \cap T_{3,2}$
equal to $1$?
\end{question}

The condition $p \geq 5$ is included here since ${\cal H}_3 \cap T_{3,2}$ is empty when $p=3$.  
When $p=5$, the authors of \cite{EP:g3small} show that ${\cal H}_3 \cap T_{3,2}$ has a unique irreducible component of dimension $2$, 
and that the generic point of this component has $p$-rank $1$.  
Thus, no component of $V_{3,0} \cap {\cal H}_3$ can be contained in $T_{g,2}$.   
As a result, the authors prove 
Conjecture \ref{Cgeneric} for $V_{g,g-3} \cap {\cal H}_g$ when $g \geq 3$ and $p=5$.

\subsection{Hyperelliptic curves when $p=2$} \label{Sp=2}

The next results show that the hypothesis $p \geq 3$ in necessary for the hyperelliptic case of Conjecture \ref{Cgeneric}.
A version of Lemma \ref{Lp=2} is stated without proof in \cite[3.2]{V:cycles}.

\begin{lemma} \label{Lp=2}
Let $p=2$.  Let $X \in V_{g,0} \cap {\cal H}_g$.  
Then $a_X=\lfloor (g+1)/2 \rfloor$. 
\end{lemma}

\begin{proof}
If $X \in V_{g,0} \cap {\cal H}_g$, then $X$ has an Artin-Schreier 
equation $y^2-y=f(x)$ for some $f(x) \in k[x]$ with ${\rm deg}(f)=2g+1$.  
A basis for $H^0(\Omega_1)$ is $\{x^bdx \ | \ 0 \leq b < g\}$.
The action of $V$ on $J_X[p]$ is the same as the action of the Cartier operator $C$ on $H^0(\Omega_1)$.
By definition (see, for example, \cite[Prop.\ 2.1]{Y}) 
the Cartier operator acts by $C(x^bdx)=0$ if $b$ is even and $C(x^bdx)=x^{(b-1)/2}dx$ if $b$ is odd.
So $\nu_g={\rm dim}(CH^0_{\Omega_1})=\#\{b \equiv 1 \bmod 2 \ | \  0 \leq b < g\}=\lceil (g-1)/2 \rceil$.
Then $a_X=g-\nu_g=\lfloor (g+1)/2 \rfloor$. 
\end{proof}

The next example shows that Conjecture \ref{Cgeneric} is false for ${\cal H}_3 \cap V_{3,1}$ when $p=2$.

\begin{example}
Let $p=2$, $g=3$, and $f=1$.  
Then ${\cal H}_{3} \cap V_{3,1}$ has two components $\Gamma_1$ and $\Gamma_2$, each of dimension $3$.
Each curve $X$ in $\Gamma_1$ has an Artin-Schreier equation of the form $y^2-y=x^5+c_1x^3+c_2x+c_3/x$.
A basis for $H^0(\Omega_1, X)$ is $\{dx/x, dx, xdx\}$.
The Cartier operator acts by $C(dx/x)=dx/x$, $C(dx)=0$ and $C(xdx)=dx$.  Thus $a_X=1$ if $X \in \Gamma_1$.
Each curve $X$ in $\Gamma_2$ has an Artin-Schreier equation of the form $y^2-y=x^3+c_1x+c_2/x+c_3/x^3$.
A basis for $H^0(\Omega_1, X)$ is $\{dx/x^2, dx/x, dx\}$.
The Cartier operator acts by $C(dx/x^2)=0$, $C(dx/x)=dx/x$ and $C(dx)=0$.  Thus $a_X=2$ if $X \in \Gamma_2$.
So Conjecture \ref{Cgeneric} is false for the component $\Gamma_2$.
\end{example}

\begin{remark}
In \cite{NS:hyper}, for $k={\mathbb F}_{2^a}$, the authors calculate the number of 
$k$-isomorphism classes of hyperelliptic curves of genus $3$ with a given Newton polygon.
\end{remark}

\section*{Acknowledgments}
The author was partially supported by NSF grant DMS-07-01303 and
would also like to thank J. Achter, D. Glass, F. Oort, and the referee.

\bibliographystyle{plain}
\bibliography{largeprank}

\end{document}